\documentclass[12pt]{article}

\usepackage{amsmath,amssymb,amsbsy}
\usepackage{graphicx}
\usepackage{diagbox,multirow}
\newtheorem{thm}{Theorem}
\newtheorem{lem}[thm]{Lemma}

\newtheorem{con}[thm]{Conjecture}
\newtheorem{obs}[thm]{Observation}

\def\gr{\mathcal{G}}

\def\gcd{\mathop{\rm gcd}\nolimits}
\newcommand{\prf}{\noindent{\it Proof.}\ }
\newcommand{\qed}{\hfill \rule{.1in}{.1in}}

\def\zet{\mathop\mathbb{Z}\nolimits}

\makeatletter
\def\imod#1{\allowbreak\mkern10mu({\operator@font mod}\,\,#1)}

\usepackage{fancyhdr}
\begin{document}

\title{Zero sum partition into sets of the same
order and its applications}

\author{Sylwia Cichacz\\
{\small Faculty of Applied Mathematics}\\
{\small AGH University of Science and Technology}\\
{\small Al. Mickiewicza 30, 30-059 Krak\'ow, Poland}\\
}
\maketitle
\begin{abstract}

We will say that an Abelian group $\Gamma$ of order $n$ has the $m$-\emph{zero-sum-partition property} ($m$-\textit{ZSP-property}) if  $m$ divides $n$, $m\geq 2$ and there is a partition of $\Gamma$ into pairwise disjoint subsets $A_1, A_2,\ldots , A_t$, such that $|A_i| = m$ and $\sum_{a\in A_i}a = g_0$ for $1 \leq i \leq t$, where $g_0$ is the identity element of $\Gamma$. \\

In this paper we study the $m$-ZSP property of $\Gamma$. We show that $\Gamma$ has $m$-ZSP if and only if $|\Gamma|$ is odd or $m\geq 3$ and $\Gamma$ has more than one involution. We will apply the results to the study of group distance magic graphs as well as to generalized Kotzig arrays.


\end{abstract}


\section{Introduction}
Assume $\Gamma$ is an Abelian group of order $n$ with the operation denoted by $+$.  For convenience
we will write $ka$ to denote $a + a + \ldots + a$ (where the element $a$ appears $k$ times), $-a$ to denote the inverse of $a$ and
we will use $a - b$ instead of $a+(-b)$. Let the identity element of $\Gamma$ be denoted by $g_0$. Recall that any group element $\iota\in\Gamma$ of order 2 (i.e., $\iota\neq0$ such that $2\iota=0$) is called an \emph{involution}.

In \cite{KLR} Kaplan, Lev and  Roditty introduced a notion of zero-sum partitions of subsets in Abelian groups. Let $\Gamma$ be an Abelian group and let $A$ be a finite subset of $\Gamma - \{g_0\}$, with $|A| = n-1$. We shall say that $A$ has the
\emph{zero-sum-partition property} (\textit{ZSP-property}) if every partition $n-1 = r_1 + r_2 + \ldots + r_t$ of $n-1$, with $r_i \geq 2$ for $1 \leq i \leq t$ and for any possible positive integer $t$, there is a partition of $A$ into pairwise disjoint subsets $A_1, A_2,\ldots , A_t$, such that $|A_i| = r_i$ and $\sum_{a\in A_i}a = g_0$ for $1 \leq i \leq t$. In the case that $\Gamma$ is finite, we shall say that $\Gamma$ has the ZSP-property if $A = \Gamma - \{g_0\}$ has the ZSP-property.

The following theorem for cyclic groups of odd order is proved.
\begin{thm}[\cite{KLR}] The group $\zet_n$ has the ZSP-property if and only if $n$ is odd.
\end{thm}
 Moreover, Kaplan, Lev and  Roditty showed that if  $\Gamma$ is a finite Abelian group of even order $n$ such that the number of involutions in $\Gamma$ is different from 3, then $\Gamma$ does not have the ZSP-property \cite{KLR}. Their results along with results proved by  Zeng \cite{Zeng}
 give necessary and sufficient conditions for the ZSP-property for  a finite Abelian group.
\begin{thm}[\cite{KLR, Zeng}] Let $\Gamma$ be a finite Abelian group. Then $\Gamma$ has the ZSP-property
if and only if either $\Gamma$ is of odd order or $\Gamma$ contains exactly three involutions. \label{ZSP}
\end{thm}
They apply those results to the study of anti-magic trees \cite{KLR, Zeng}.\\

Consider a simple graph $G$ whose order we denote by $n=|G|$. We denote by $V (G)$  the vertex set and $E(G)$  the edge set of a graph $G$. The \emph{open neighborhood} $N(x)$ of a vertex $x$ is the set of vertices adjacent to $x$, and the degree $d(x)$ of $x$ is $|N(x)|$, the size of the neighborhood of $x$. \\

In this paper we investigate also distance magic labelings, which belong to a large family of magic-type labelings. Generally speaking, a magic-type labeling of a graph $G=(V,E)$ is a mapping from $V, E,$ or $V\cup E$ to a set of labels which most often is a set of integers or group elements. Then the weight of a graph element is typically the sum of labels of the neighboring elements of one or both types. If the weight of each element is required to be equal, then we speak about magic-type labeling; when the weights are all different (or even form an arithmetic progression), then we speak about an antimagic-type labeling. Probably the best known problem in this area is the {\em antimagic conjecture} by Hartsfield and Ringel~\cite{HR}, which claims that the edges of every graph except $K_2$ can be labeled by integers $1,2,\dots,|E|$ so that the weight of each vertex is different.

A \emph{distance magic labeling} (also called \emph{sigma
labeling}) of a graph $G=(V,E)$ of order $n$ is a bijection $\ell
\colon V \rightarrow \{1, 2,\ldots  , n\}$ with the property that
there is a positive integer $k$ (called the \emph{magic constant}) such
that
$$w(x)=\sum_{y\in N(x)}\ell(y)=k~ {\text{for every}}~ x \in V(G),$$
where $w(x)$ is the {\em weight} of vertex $x$. If a graph $G$ admits a distance magic labeling, then we say that $G$ is a \emph{distance magic graph}.

The following result was proved in \cite{MRS}.
\begin{obs}[\cite{MRS}] There is no distance magic $r$-regular graph with $r$ odd.\label{nieparzyste}\end{obs}

Froncek in \cite{Fro} defined the notion of group distance magic graphs, i.e., the graphs allowing a bijective labeling of vertices
with elements of an Abelian group resulting in constant sums of neighbor labels.\\

A $\Gamma$\emph{-distance magic labeling} of a graph $G = (V, E)$ with $|V| = n$ is a bijection $\ell$ from $V$ to an Abelian group $\Gamma$ of order $n$
such that the weight $w(x) =\sum_{y\in N(x)}\ell(y)$ of every vertex $x \in V$ is equal to the same element $\mu\in \Gamma$, called the \emph{magic
constant}. A graph $G$ is called a \emph{group distance magic graph} if there exists a $\Gamma$-distance magic labeling for every Abelian
group $\Gamma$ of order $|V(G)|$.\\

Notice that the constant sum partitions of a group $\Gamma$ lead to complete multipartite $\Gamma$-distance magic labeled graphs. For
instance, the partition $\{0\}$, $\{1, 2, 4\}$, $\{3, 5, 6\}$ of the group $\zet_7$ with constant sum $0$ leads to a $\zet_7$-distance magic labeling
of the complete tripartite graph $K_{1,3,3}$ (see \cite{Cic2,Cic3}).

The connection between distance magic graphs and $\Gamma$-distance magic graphs is as follows. Let $G$ be a distance magic
graph of order $n$ with the magic constant $\mu'$. If we replace the label $n$ in a distance magic labeling for the graph $G$
by the label $0$, then we obtain a $\zet_n$-distance magic labeling for the graph $G$ with the magic constant $\mu = \mu'\pmod n$.
Hence every distance magic graph with $n$ vertices admits a $\zet_n$-distance magic labeling. However a $\zet_n$-distance magic
graph on $n$ vertices is not necessarily a distance magic graph. Moreover, there are some
graphs that are not distance magic while at the same time they are group distance magic (see \cite{Cic2}).

A general theorem for $\Gamma$-distance magic labeling similar to Observation~\ref{nieparzyste} was proved recently.
\begin{thm}[\cite{CicFro}]\label{gr:odd} Let $G$ be an $r$-regular graph on $n$
vertices, where $r$ is odd.
There does not exist an Abelian group $\Gamma$ of order $n$ with exactly one
involution $\iota$ such that $G$ is $\Gamma$-distance magic.
\end{thm}

The paper is organized as follows. In Section~\ref{sZSP}  we show that $\Gamma$ has $m$-ZSP if and only if $|\Gamma|$ is odd or $m\geq 3$ and $\Gamma$ has more than one involution.
So far there is known only one example (namely $K_{3,3,3,3}$, \cite{CicFro}) having a $\Gamma$-distance magic labeling for a group $\Gamma\cong\zet_2\times\zet_2\times\zet_3$ which has more than one involution. In Section~\ref{GDM} using $m$-ZSP of Abelian groups we show an infinite family of odd regular graphs possessing $\Gamma$-distance magic labeling for groups $\Gamma$ that have more than one involution.  Finally in Section~ \ref{Kotzig} we introduce a generalization of Kotzig array and we give necessary and sufficient conditions for their existence.

The following theorem was proved in \cite{Cic3}.
\begin{thm}[\cite{Cic3}]\label{2mod4_1} Let $G$ have order $n \equiv 2 \pmod 4$
with all vertices having odd degree.
There does not exist an Abelian group $\Gamma$ of order $n$ such that $G$ is a $\Gamma$-distance magic graph.
\end{thm}

\section{Preliminaries}

 A non-trivial
finite group has elements of order $2$ (an involution) if and only if the order of the group is even. The fundamental theorem of finite Abelian groups states that a finite Abelian
group $\Gamma$ of order $n$ can be expressed as the direct product of cyclic subgroups of prime-power order. This implies that
$$\Gamma\cong\zet_{p_1^{\alpha_1}}\times\zet_{p_2^{\alpha_2}}\times\ldots\times\zet_{p_k^{\alpha_k}}\;\;\; \mathrm{where}\;\;\; n = p_1^{\alpha_1}\cdot p_2^{\alpha_2}\cdot\ldots\cdot p_k^{\alpha_k}$$
and $p_i$ for $i \in \{1, 2,\ldots,k\}$ are not necessarily distinct primes. This product is unique up to the order of the direct product. When $t$ is the number of these cyclic components
whose order is a multiple of $2$, then $\Gamma$ has $2^t-1$ involutions. In particular, if $n \equiv 2 \pmod 4$, then $\Gamma\cong \zet_2\times \Lambda$ for some
Abelian group $\Lambda$ of odd order $n/2$.  Moreover every cyclic group of even order has exactly one involution. Since the properties and results in this paper are invariant under the isomorphism between groups, we only need to consider a group in an isomorphism class. 
The sum of all the
group elements is equal to the sum of the involutions and the neutral element.

The following lemma was proved in~\cite{CN} (see \cite{CN}, Lemma 8).

\begin{lem}[\cite{CN}]\label{involutions} Let $\Gamma$  be an Abelian group.
\begin{itemize}
 \item[-] If $\Gamma$ has exactly one involution $\iota$, then $\sum_{g\in \Gamma}g= \iota$.
\item[-] If $\Gamma$ has no involutions, or more than one involution, then $\sum_{g\in \Gamma}g=g_0$.
\end{itemize}
\end{lem}

Zeng proves a lemma which plays an important role in the proof of the
main result  (see \cite{Zeng}, Lemma 2.1).

\begin{lem}[\cite{Zeng}]\label{bijection}
Let $\Gamma$ be a finite Abelian group of odd order or $\Gamma$ contains exactly three involutions and Bij$(\Gamma)$ denote the set of all bijections from $\Gamma$ to itself.
Then there exist $\phi,\varphi\in$Bij$(\Gamma)$ (not necessarily distinct) such that $g+\phi(g)+\varphi(g)=g_0$  for every $g\in\Gamma$. 
\end{lem}

A distance magic graph $G$ with an even number of vertices is called \textit{balanced} if
there exists a bijection $\ell \colon V(G)\rightarrow\{1,\ldots,|V(G)|\}$ such that for every $w\in V(G)$ the following holds: if $u\in N(w)$ with $\ell(u)=i$, then there exists $v\in N(w)$, $v\neq u$,
with $\ell(v)=|V(G)|+1-i$.   It is easy to see that a balanced distance magic graph has an even number
of vertices and that it is an $r$-regular graph for an even $r$. Simple examples are empty graph
on an even number of vertices $\overline{K}_{2n}$, cycle $C_4$, and $K_{2n}-I$, for a perfect matching $I$. 

The \textit{lexicographic product} or \textit{graph composition} $G \circ H$ of graphs $G$ and $H$ is a graph such that  the vertex set of $G \circ H$ is the Cartesian product $V(G) \times V(H)$; and
    any two vertices $(u,v)$ and $(x,y)$ are adjacent in $G \circ H$ if and only if either $u$ is adjacent with $x$ in $G$ or $u = x$ and $v$ is adjacent with $y$ in $H$. Note that $G\circ H$ and $H\circ G$ are not isomorphic in general. One can imagine obtaining $G\circ H$ by blowing up each vertex of $G$ into a copy of $H$.

It is interesting that if you replace each vertex of an
$r$-regular graph by some specific $p$-regular graph (\cite{ref_AnhCicPetTep2,ref_AnhCicPetTep1}), then the obtained graph is a group distance magic graph. The
following theorem was proved in \cite{ref_AnhCicPetTep2} and will be used in the fourth section.

\begin{thm}[\cite{ref_AnhCicPetTep2}]Let $G$ be a regular graph and $H$ be a graph not isomorphic to $\overline{K}_n$ where
$n$ is odd. Then $G\circ H$ is a balanced distance magic graph if and only if $H$ is a balanced
distance magic graph.\label{tw1}
\end{thm}

\begin{thm}[\cite{ref_AnhCicPetTep1}] If $G$ is a balanced distance magic graph, then $G$ is a group distance magic
graph.\label{tw2}
\end{thm}

\section{Zero sum partition}\label{sZSP}
Let the identity element of $\Lambda$ be denoted by $a_0$. By $\langle g \rangle$ we denote o subgroup generated by $g$ in the group $\Gamma$. A quotient group $\Gamma$ modulo $H$ for a subgroup $H$ of $\Gamma$ will be denoted by $\Gamma/N$.\\

For convenience, let $\gr$ denote the set consisting of all Abelian groups which are
of odd order or contain more than one involution.

We will say that an Abelian group $\Gamma$ of order $n$ has the $m$-\emph{zero-sum-partition property} ($m$-\textit{ZSP-property}) if  $m$ divides $n$, $m\geq 2$ and there is a partition of $\Gamma$ into pairwise disjoint subsets $A_1, A_2,\ldots , A_t$, such that $|A_i| = m$ and $\sum_{a\in A_i}a = g_0$ for $1 \leq i \leq t$. \\
We will start with some lemmas.
\begin{lem}\label{involutions3} Let $\Gamma$  be an Abelian group  with involutions set $I^* =\{\iota_1,\iota_2,\ldots,\iota_{2^p-1}\}$, $p >1$ and let $I = I^*\cup\{g_0\}$.
 There exists a partition $A=\{A_1, A_2,\ldots, A_{2^{p-2}}\}$ of $I$ such that $|A_i|=4$,
   $\sum_{a\in A_i}a=g_0$ for $i\in\{1,2,\ldots,2^{p-2}\}$.
\end{lem}
\prf Let $\iota_0=g_0$. Recall that since $I_p=\{\iota_0,\iota_1,\ldots,\iota_{2^p-1}\}$ is a subgroup of $\Gamma$, we have $I_p\cong (\zet_2)^p$. One can check that the lemma is true for $p=2$. The sufficiency will be proved then  by induction on $2$. Namely, suppose the assertion is true for some $p\geq 2$. We want to prove it is true for
$r=p+ 1$.  By the inductive hypothesis the group $I_{p-1}\cong (\zet_2)^{p-1}$ has $4$-ZSP. Let  $A_1,A_2,\ldots,A_{2^{p-3}}$ be the desired partition of  $I_{p-1}$. Define $A_i^0$ by replacing each element $a$ in $\Lambda$ by $(0,a)\in I_p$  and analogously $A_i^1$ by replacing each element $a$ in $\Lambda$ by $(1,a)\in I_p$. Obviously $\sum_{a\in A_i^j}a=g_0$ for $j=0,1$, $i=1,2,\ldots,2^{p-3}$~\qed\\

Using the same method as in proof of Lemma \ref{bijection}, we can obtain the following lemma.
\begin{lem}\label{bijections}
Let $\Gamma$ be a finite Abelian group of odd order or $\Gamma$ contains more than one involution and Bij$(\Gamma)$ denote the set of all bijections from $\Gamma$ to itself.
Then there exist $\phi,\varphi\in$Bij$(\Gamma)$ (not necessarily distinct) such that $g+\phi(g)+\varphi(g)=g_0$  for every $g\in\Gamma$. 
\end{lem}
\prf First we prove an assertion: let $\Gamma_1,\Gamma_2\in\gr$ and suppose we have proven
the lemma for $\Gamma_1,\Gamma_2$, then the lemma holds for $\Gamma=\Gamma_1\times\Gamma_2$. Indeed, suppose there bijections for $\Gamma_1$ are $\phi_1$ and $\varphi_1$, whereas for $\Gamma_2$ $\phi_2$ and $\varphi_2$. Then
$$\phi=(\phi_1,\phi_2): \Gamma_1\times \Gamma_2\rightarrow\Gamma_1\times \Gamma_2,\;\; (a_1,a_2)\mapsto (\phi_1(a_1),\phi_2(a_2))$$
and
$$\varphi=(\varphi_1,\varphi_2): \Gamma_1\times \Gamma_2\rightarrow\Gamma_1\times \Gamma_2,\;\; (a_1,a_2)\mapsto (\varphi_1(a_1),\varphi_2(a_2)).$$
By the assertion above and Lemma~\ref{bijection} and noting that this lemma is invariant under the isomorphism, it suffices to prove the lemma for $\Gamma=\zet_{2^\alpha}\times\zet_{2^\beta}\times\zet_{2^\kappa}$  with $\alpha,\beta,\kappa\geq1$.\\
The proof is by induction on $|G|$. We handle three basic cases.  For $\zet_{2}\times\zet_{2}\times\zet_{2}$   the table below  gives the desired bijections.\\
\begin{tabular}{|c|c|c|c|c|c|c|c|c|}
\hline
$g$          & (0,0,0) & (0,0,1) & (0,1,0) & (0,1,1) & (1,0,0) & (1,0,1) & (1,1,0) & (1,1,1)\\ \hline
$\phi(g)$    & (0,0,1) & (1,1,1) & (0,1,0) & (1,0,0) & (0,0,0) & (1,1,0) & (0,1,1) & (1,0,1)\\ \hline
$\varphi(g)$ & (0,0,1) & (1,1,0) & (0,0,0) & (1,1,1) & (1,0,0) & (0,1,1) & (1,0,1) & (0,1,0)\\ \hline
\end{tabular}\newline

Whereas  for  $\zet_{2}\times\zet_{2}\times\zet_{4}$  the table:\\
\begin{tabular}{|c|c|c|c|c|c|c|c|c|}
\hline
$g$          & (0,0,0) & (0,1,0) & (1,0,0) & (1,1,0) & (0,0,1) & (0,1,1) & (1,0,1) & (1,1,1)\\ \hline
$\phi(g)$    & (0,1,0) & (1,1,2) & (1,0,1) & (0,0,3) & (1,0,3) & (0,0,1) & (0,1,3) & (1,1,1)\\ \hline
$\varphi(g)$ & (0,1,0) & (1,0,2) & (0,0,3) & (1,1,1) & (1,0,0) & (0,1,2) & (1,1,0) & (0,0,2)\\\hline\hline
$g$          & (0,0,2) & (0,1,2) & (1,0,2) & (1,1,2) & (0,0,3) & (0,1,3) & (1,0,3) & (1,1,3)\\ \hline
$\phi(g)$    & (0,0,2) & (1,0,0) & (1,1,3) & (0,1,1) & (0,0,0) & (1,0,2) & (1,1,0) & (0,1,2)\\ \hline
$\varphi(g)$ & (0,0,0) & (1,1,2) & (0,1,3) & (1,0,1) & (0,0,1) & (1,1,3) & (0,1,1) & (1,0,3)\\ \hline
\end{tabular}
\newline\newline

Suppose that $\zet_{2}\times\zet_{2}\times\zet_{8}$ , for $g=(i,j,0)$ and $g=(i,j,4)$ we have the following table.\\
\begin{tabular}{|c|c|c|c|c|c|c|c|c|}
\hline
$g$          & (0,0,0) & (0,1,0) & (1,0,0) & (1,1,0) & (0,0,4) & (0,1,4) & (1,0,4) & (1,1,4)\\ \hline
$\phi(g)$    & (0,1,0) & (1,1,4) & (1,0,0) & (0,0,4) & (0,0,0) & (1,0,4) & (1,1,0) & (0,1,4)\\ \hline
$\varphi(g)$ & (0,1,0) & (1,0,4) & (0,0,0) & (1,1,4) & (0,0,4) & (1,1,4) & (0,1,4) & (1,0,0)\\ 
\hline
\end{tabular}\newline

If $g=(i,j,l)$ for $l \not \in\{0,4\}$, then set a triple $(a,b,c)\in \{(2,3,3),(7,2,7),$ $(5,6,5),(6,1,1),(1,5,2),(3,7,6)\}$, and for $g=(a,i,j)$ we have the following table.\\
\begin{tabular}{|c|c|c|c|c|}
\hline$g$    & (0,0,a) & (0,1,a) & (1,0,a) & (1,1,a) \\ \hline
$\phi(g)$    & (0,0,b) & (1,0,b) & (1,1,b) & (0,1,b) \\ \hline
$\varphi(g)$ & (0,0,c) & (1,1,c) & (0,1,c) & (1,0,c) \\ \hline

\end{tabular}\newline

Now we proceed to the induction part. 
Suppose first that $\Gamma=\zet_{2^\alpha}\times\zet_{2^\beta}\times\zet_{2^\kappa}$  with $\beta,\kappa\geq2$. Then there exists a subgroup $\Gamma_0=\langle 1 \rangle\times \langle 2\rangle\times \langle 2 \rangle\in\gr$ of $\Gamma$ such that $\Gamma/\Gamma_0\cong \zet_2\times\zet_2$.
By the induction hypothesis, there are $\phi_0$, $\varphi_0\in$Bij$(\Gamma_0)$ such that $a+\phi_0(a)+\varphi_0(a)=g_0$ for every $a\in G_0$. Choose a set of coset representatives for $\Gamma/\Gamma_0$ to be $\{0,c,d,-c-d\}$. Note that
$$(c+\Gamma_0)\bigcup(d+\Gamma_0)\bigcup(-c-d+\Gamma_0)=\bigcup_{b\in\Gamma_0}\{c+b,d+\phi_0(b),-c-d+\varphi_0(b)\},$$
and every subset $\{c+b,d+\phi_0(b),-c-d+\varphi_0(b)\}$ is zero-sum. Thus define $\phi$ and
$\varphi$ as follows.
$$\phi(a)=\left\{\begin{array}{ccl}
                     \phi_0(a) & \mathrm{if} & a\in\Gamma_0, \\
											d+\phi_0(b) & \mathrm{if} & a=c+b \;\;\mathrm{for}\;\; \rm{some}\;\; b\in \Gamma_0, \\
											-c-d+\varphi_0(b) & \mathrm{if} & a=d+\phi_0(b)\;\;\mathrm{for}\;\; \rm{some}\;\; b\in \Gamma_0, \\
												c+b & \mathrm{if} & a=-c-d+\varphi_0(b) \;\;\mathrm{for}\;\; \rm{some}\;\; b\in \Gamma_0, \\
 \end{array}\right.$$
and
$$\varphi(a)=\left\{\begin{array}{ccl}
                     \varphi_0(a) & \mathrm{if} & a\in\Gamma_0, \\
											-c-d+\varphi_0(b) & \mathrm{if} & a=c+b \;\;\mathrm{for}\;\; \rm{some}\;\; b\in \Gamma_0, \\
											c+b & \mathrm{if} & a=d+\phi_0(b) \;\;\mathrm{for}\;\; \rm{some}\;\; b\in \Gamma_0, \\
												d+\phi_0(b) & \mathrm{if} & a=-c-d+\varphi_0(b) \;\;\mathrm{for}\;\; \rm{some}\;\; b\in \Gamma_0. \\
 \end{array}\right.$$
Suppose now that $\Gamma=\zet_{2}\times\zet_{2}\times\zet_{2^\kappa}$. The cases when $\kappa \leq 3$ have been shown in the base of the induction. Thus we may assume that  $\kappa\geq 4$ and there exists a subgroup $\Gamma_0=\langle 1 \rangle\times \langle 1 \rangle\times \langle 8 \rangle\in\gr$ such that $\Gamma/\Gamma_0\cong \zet_8$. By the induction hypothesis, there are $\phi_0$, $\varphi_0\in$Bij$(\Gamma_0)$ such that $a+\phi_0(a)+\varphi_0(a)=g_0$ for every $a\in \Gamma_0$. Since
$$\zet_8=\{0,4\}\bigcup\{1,2,5\}\bigcup\{-1,-2,-5\},$$
we can choose a set of coset representatives for $\Gamma/\Gamma_0$ to be $\{0,e,c,d,-c-d,-c,-d,c+ d\}$ with $2e\in \Gamma_0$. Since $\Gamma_1=\Gamma_0\cup (e+\Gamma_0)$ (for example $e=(0,0,4)$ and then $\Gamma_1=\langle 1 \rangle\times \langle 1 \rangle\times \langle 4 \rangle$) is a group and moreover $\Gamma_0$ is a subgroup of $\Gamma_1$ and $\Gamma_1$ is a subgroup of $\Gamma$, we have $\Gamma_1\in\gr$. Thus there exist $\phi_1,\varphi_1\in$Bij$(\Gamma_1)$ such that $a+\phi_1(a)+\varphi_1(a)=g_0$ for every $a\in \Gamma_1$. Similarly as above, we have
$$(c+\Gamma_0)\bigcup(d+\Gamma_0)\bigcup(-c-d+\Gamma_0)=\bigcup_{b\in\Gamma_0}\{c+b,d+\phi_0(b),-c-d+\varphi_0(b)\},$$
and
$$(-c+\Gamma_0)\bigcup(-d+\Gamma_0)\bigcup(c+d+\Gamma_0)=\bigcup_{b\in\Gamma_0}\{-c+b,-d+\phi_0(b),c+d+\varphi_0(b)\}.$$
Hence $\Gamma=\Gamma_1\cup(\cup_{i=1}^{2|\Gamma_0|}\{a_{i,1},a_{i,2},a_{i,3}\})$ where  $a_{i,1}+a_{i,2}+a_{i,3}=g_0$ for $1\leq i\leq 2|\Gamma_0|$. Define $\phi$ and
$\varphi$ as follows.
$$\phi(a)=\left\{\begin{array}{ccl}
                     \phi_1(a) & \mathrm{if} & a\in\Gamma_1, \\
											a_{i,2} & \mathrm{if}& a=a_{i,1}\;\;\mathrm{for}\;\; \rm{some}\;\; 1\leq i\leq2|\Gamma_0|, \\
											a_{i,3} & \mathrm{if}& a=a_{i,2} \;\;\mathrm{for}\;\; \rm{some}\;\;1\leq i\leq2|\Gamma_0|, \\
											a_{i,1} & \mathrm{if} &a=a_{i,3}\;\;\mathrm{for}\;\; \rm{some}\;\; 1\leq i\leq2|\Gamma_0|, \\
 \end{array}\right.$$
and
$$\varphi(a)=\left\{\begin{array}{ccl}
                    \varphi_1(a) & \mathrm{if} & a\in\Gamma_1, \\
											a_{i,3} & \mathrm{if} &a=a_{i,1}\;\;\mathrm{for}\;\; \rm{some}\;\; 1\leq i\leq2|\Gamma_0|, \\
											a_{i,1} & \mathrm{if} &a=a_{i,2} \;\;\mathrm{for}\;\; \rm{some}\;\; 1\leq i\leq2|\Gamma_0|, \\
											a_{i,2} & \mathrm{if}& a=a_{i,3} \;\;\mathrm{for}\;\; \rm{some}\;\; 1\leq i\leq2|\Gamma_0|. \\
 \end{array}\right.$$
This completes the proof of the lemma.~\qed

\begin{lem}\label{mZSPeven} An Abelian group $\Gamma$  of order $n$ has $(2m)$-ZSP-property  if and only if $(2m)|n$, $m\geq2$ and $\Gamma\in\gr$.
\end{lem}
\prf We will show the necessity first. Suppose that $\Gamma$  has $(2m)$-ZSP-property and $A_1, A_2,\ldots , A_t$ is the desired partition. Therefore the order of $\Gamma$ is even divisible by $2m$. Hence it has at least one involution. Assume now, that there is exactly one involution $\iota$. On one hand $\sum_{g\in \Gamma}g=\sum_{i=1}^t\sum_{a\in A_i}a=g_0$, on the other $\sum_{g\in \Gamma}g=\iota$ by Lemma~\ref{involutions}, a contradiction.

Let $\Gamma$  be an Abelian group  with involutions set $I^* =\{\iota_1,\iota_2,\ldots,\iota_{2^p-1}\}$, $p >1$  and let $I = I^*\cup\{g_0\}$. There exists a partition $\{I_1, I_2,\ldots, I_{2^{p-2}}\}$ of $I$ such that $|I_i|=4$,
   $\sum_{a\in I_i}a=g_0$ for $i\in\{1,2,\ldots,2^{p-2}\}$. 
Recall that $\Gamma$ can be expressed as a direct product of cyclic subgroups of prime-power order.  When $p$ is the number of these cyclic components
whose order is a power of $2$, then $\Gamma$ has $2^p-1$ involutions.	
	Note that all element in the set $\Gamma\setminus I=\{g_1,g_2,\ldots,g_{|\Gamma|-2^p}\}$ can be partition in the pairs $B_i=\{g_i,-g_i\}$ for $i=1,2,\ldots,(|\Gamma|-2^p)/2$.

It is easily seen now that for $2m\equiv 0 \pmod 4$ there exists a partition $A=\{A_1, A_2,\ldots, A_{|\Gamma|/2m}\}$ of $\Gamma$ such that $|A_i|=2m$,
   $\sum_{a\in A_i}a=g_0$ for $i\in\{1,2,\ldots,|\Gamma|/2m\}$. 

Note that for $2m\equiv 2 \pmod 4$, because $m\geq3$ is odd $\Gamma\cong\zet_{2^{\alpha_1}}\times\zet_{2^{\alpha_2}}\times\ldots\times\zet_{2^{\alpha_p}}\times\Delta$ for $\alpha_1,\alpha_2,\ldots,\alpha_p\geq1$ and some Abelian group $\Delta$ of odd order divisible by $m$. Therefore $(|\Gamma|-2^p)=|\Gamma\setminus I|\geq (m-1)|I|=(m-1)2^p$. Hence $(|\Gamma|-2^p)/2\geq2^p>2^{p-2}$ and one can check that there exists a partition $A=\{A_1, A_2,\ldots, A_{|\Gamma|/2m}\}$ of $\Gamma$ such that $|A_i|=2m$,
   $\sum_{a\in A_i}a=g_0$ for $i\in\{1,2,\ldots,|\Gamma|/2m\}$.~\qed

\begin{lem}\label{mZSPodd} An Abelian group $\Gamma$  of order  $n$ has $(2m+1)$-ZSP-property  if and only if  $m\geq1$, $(2m+1)|n$ and $\Gamma\in\gr$.
\end{lem}
\prf We will show the necessity first. The assumption  $m\geq1$ is obvious. Suppose that $\Gamma$  has one involution and it has $(2m+1)$-ZSP-property. Let $A_1, A_2,\ldots , A_t$ be the desired partition. Therefore the order of $\Gamma$ is even. Hence it has at least one involution. On one hand $\sum_{g\in \Gamma}g=\sum_{i=1}^t\sum_{a\in A_i}a=g_0$, on the other $\sum_{g\in \Gamma}g=\iota$ by Lemma~\ref{involutions}, contradiction.

If $|\Gamma|$ is odd we are done by Theorem~\ref{ZSP}. Suppose now that $\Gamma$ has  $2^{p}-1$ involutions for some positive integer $p$. Therefore $\Gamma\cong\Gamma_0\times\Lambda$ for  $\Gamma_0\cong\zet_{2^{\alpha_1}}\times\zet_{2^{\alpha_2}}\times\ldots\times\zet_{2^{\alpha_p}}$ for $\alpha_1,\alpha_2,\ldots,\alpha_p\geq1$ and some Abelian group $\Delta$  of odd order divisible by $(2m+1)$.

 Hence $\Lambda$ has $(2m+1)$-ZSP-property and there  exists a partition $A_1,A_2,\ldots,$ $A_{t}$ of the group $\Lambda$ such that for every $1 \leq i \leq t$, $|A_i| = 2m+1$ with $\sum_{a\in A_i} a =a_0$. Denote the $h$-th element in the set $A_i$ by $a_{i,h}$.\\
Let $\gamma(w)=(2n-1)w$ for $w\in\Gamma_0$. Note that  $\gamma\in$Bij$(\Gamma_0)$, because $\gcd(2^{\alpha},2m-1)=1$ for any positive integer $\alpha$. Recall that there exist $\phi,\varphi\in$Bij$(\Gamma_0)$  such that $g+\phi(g)+\varphi(g)=g_0^0$  for every $g\in\Gamma_0$ ($g_0^0$ is the identity element of $\Gamma_0$) by Lemma \ref{bijection}.  Let  $f\colon\Gamma_0\rightarrow\zet_{|\Gamma_0|}$ be a bijection. Define $A_i^j$ by replacing an element $a_{i,h}$ in $\Lambda$ by $a_{i,1}^j=(\phi(f(j)),a_{i,1})$, $a_{i,2}^j=(\varphi(f(j)),a_{i,2})$, $a_{i,h}^j=(\gamma^{-1}(f(j)),a_{i,h})$ for $j=0,1,\ldots,|\Gamma_0|-1$, $i=1,2,\ldots,t$, $h=3,4,\ldots,2m+1$. Observe that $\sum_{h=3}^{2m-1}a_{i,h}^j=(f(j),a_0)$ for $j=0,1,\ldots,|\Gamma_0|-1$, $i=1,2,\ldots,t$.  Therefore  $A_1^0,\ldots,A_t^0,A_1^1,\ldots,A_t^1,\ldots,A_1^{|\Gamma_0|-1},$ $\ldots,A_t^{|\Gamma_0|-1}$ is a partition of the group $\Gamma$ such that for every $1 \leq i \leq t$, $0\leq j\leq|\Gamma_0|-1$ $|A_i^j| = 2m+1$ with $\sum_{a\in A_i^j} a =(g_0^0,a_0)=g_0$.~\qed\\

By above Lemmas~\ref{mZSPeven} and \ref{mZSPodd} the following theorem is true.
\begin{thm}\label{mZSP} An Abelian group $\Gamma$  of order even $n$  has $m$-ZSP-property  if and only if  $m\geq3$, $m|n$ and $\Gamma\in\gr$.~\qed
\end{thm}

\section{Group distance labeling of odd regular graphs}\label{GDM}

\begin{thm}\label{kn} If $G$ is a graph of order $t$, then the lexicographic product $G\circ \overline{K}_{n}$ has a $\Gamma$-distance magic labeling for $n\geq3$ and any $\Gamma\in\gr$ of order $nt$.
\end{thm}
\prf Let $G$ be a  graph with the vertex set $V(G) = \{v_0,v_1, \ldots, v_{t-1} \}$, whereas  $V(\overline{K}_{2n+1}) = \{x_0,x_1,\ldots, x_{n-1}\}$.  Since $\Gamma\in\gr$,  there exists a partition $A_1,A_2,\ldots,A_t$ of  $\Gamma$ such that   $|A_i| =n$ with $\sum_{a\in A_i} a =g_0$ for every $1 \leq i \leq t$ by Theorem~\ref{mZSP}.  Label the vertices of the  $i$-th copy of $\overline{K}_{n}$ using elements from the set $A_i$ for $i=1,2,\ldots,t$.\\
 Notice that $\sum_{j=1}^{n}{\ell}(v_i,x_j)=g_0$  for $i=1,2,\ldots,t$. Therefore $w(v_i,x_j)=g_0$.~\qed

\begin{obs}\label{kn} If $G$ is a graph of odd order $t$, then the lexicographic product $G\circ \overline{K}_{2n+1}$ is  group distance magic for $n\geq1$.~\qed
\end{obs}

It is worthy to mention that by Theorems \ref{tw1} and \ref{tw2} we have the following easy observation.
\begin{obs}\label{Iztok} If $G$ is a regular graph, then the lexicographic product $G\circ \overline{K}_{2n}$  is group distance magic labeling for any $n\geq1$.~\qed
\end{obs}

Observe that if $G$ is an odd regular graph, then  the lexicographic product $G\circ \overline{K}_{2n+1}$ is also an odd regular graph. Thus the order of the graph $G\circ \overline{K}_{2n+1}$ is even by handshaking lemma. Recall that $K_{3,3,3,3}\cong K_4\circ\overline{K}_{3}$ which is $\zet_2\times\zet_2\times\zet_3$-distance magic. From the above Theorems \ref{kn} and \ref{gr:odd} we obtain the following theorem showing that there exist infinitely many odd regular  $\Gamma$-distance magic graphs for a group $\Gamma$ having more than one involution. 
\begin{thm}If $G$ is an odd regular graph of order $t$, then the lexicographic product $G\circ \overline{K}_{2n+1}$ has a $\Gamma$-distance magic labeling for $n\geq1$ and $|\Gamma|=(2n+1)t$ if and only if $\Gamma$ has more than one involution.~\qed
\end{thm}

For bipartite Eulerian graphs (i.e. all vertices have even degree)  such that the partition sets have the same cardinality we have the following theorem.
\begin{thm} If $G$ is an Eulerian bipartite graph of order $t\equiv 2 \pmod 4$ with partition sets of the same cardinality, then $G\circ \overline{K}_{2n+1}$  is group distance magic for $n\geq1$. \label{dwudzielne}
\end{thm}
\prf Note that  the partition sets $V_1$, $V_2$ of $G$ has the same cardinality $k=t/2$.
Let  $V_1 = \{v_0,v_1, \ldots, v_{k-1} \}$, $V_2 = \{w_0,w_1, \ldots, w_{k-1} \}$ whereas  $V(\overline{K}_{2n+1}) = \{x_0,x_1,\ldots, x_{2n}\}$.
Let $\Gamma$ be any Abelian group of order $(2n+1)t$. Observe that $\Gamma\cong \zet_2\times \Lambda$ for some Abelian group $\Lambda$ of odd order $(2n+1)t/2$ because $t\equiv 2 \pmod 4$. \\
Since $2n+1\geq3$,  there exists a partition $A_1,A_2,\ldots,A_{k}$ of the group $\Lambda$ such that for every $1 \leq i \leq t/2$, $|A_i| = 2n+1$ with $\sum_{a\in A_i} a =a_0$ by Theorem~\ref{mZSP}. Define $A_i^0$ by replacing each element $a$ in $\Lambda$ by $(0,a)$  and analogously $A_i^1$ by replacing each element $a$ in $\Lambda$ by $(1,a)$. Label the vertices of the  $i$-th copy of $\overline{K}_{2n+1}$ in $V_1$, $V_2$ using elements from the set $A_i^0$, $A_i^1$, respectively for $i=1,2,\ldots,k$.\\
 Notice that $w(x)=g_0$ for any $x\in G\circ \overline{K}_{2n+1}$ because $x$ has even degree.~\qed\\

Recall that any regular  bipartite graph $G$ has the partition sets with the same cardinality, therefore by above Observation~\ref{Iztok} and Theorem~\ref{dwudzielne} we easily obtain the following observation. 

\begin{obs} If $G$ is an $r$-regular bipartite   graph of order $t\equiv 2 \pmod 4$ and $n\geq2$, then $G\circ \overline{K}_{n}$  is group distance magic if and only if $rn$ is even.
\end{obs}
\prf Assume first that $rn$ is odd, then $nt\equiv 2 \pmod 4$ and by Theorem~\ref{gr:odd} there does not exist a group $\Gamma$ such that $G\circ \overline{K}_{n}$ admits a $\Gamma$-distance labeling.
If $n$ is even then we are done by Observation~\ref{Iztok}. In the last case is $n$ odd and $r$ even we apply Theorem~\ref{dwudzielne}.~\qed

\section{$\Gamma$-Kotzig arrays}\label{Kotzig}
In \cite{Wal} Marr and Wallis give a definition of a Kotzig array as a $j\times k$ grid, each row being a permutation of $\{0,1,\dots,m-1\}$ and each column having the same sum.
 \begin{lem}[\cite{Wal}]\label{array}
A Kotzig array of size $j \times k$ exists whenever $j>1$ and $j(k-1)$ is even.
\end{lem}
Note that a  Latin square of size $n \times n$ is a special case of Kotzig array of size $n \times n$.

Kotzig arrays play important role for graph labeling.  There are many constructions based on Kotzig arrays of various magic-type labelings of copies and products of regular graphs. 

For an Anelian group $\Gamma$ of order $k$ we define a $\Gamma$-Kotzig array  of size $j\times k$ as a  $j\times k$ grid, each row being a permutation of $\Gamma$ and each column having the same sum. Denote by $x_{i,j}$ the entry in the $i$-th row and $j$-th column of
the array.
\begin{lem}
A $\Gamma$-Kotzig array of size $2 \times k$ exists for any Abelian $\Gamma$.
\end{lem}
\prf Let $\Gamma=\{g_0,g_1,\ldots,g_{k-1}\}$. Let $x_{1,i}=g_i$, $x_{2,i}=-g_i$ for $i=0,1,\ldots,k-1$.~\qed

\begin{lem}
A $\Gamma$-Kotzig array of size $3 \times k$ exists for any Abelian $\Gamma\in\gr$.
\end{lem}
\prf Let $\Gamma=\{g_0,g_1,\ldots,g_{k-1}\}$.  There exist $\phi,\varphi\in$Bij$(\Gamma)$  such that $g+\phi(g)+\varphi(g)=g_0$  for every $g\in\Gamma$ by Lemma~\ref{bijections}. Let $x_{1,i}=g_i$, $x_{2,i}=\phi(g_i)$, $x_{2,i}=\varphi(g_i)$.~\qed

 \begin{thm}
A $\Gamma$-Kotzig array of size $j \times k$ exists whenever $j>1$ and $j$ is even or $\Gamma\in \gr$.
\end{thm}
\prf Obviously for $j=1$ there does not exist a $\Gamma$-Kotzig array. 
 Assume first $j$ is even. To construct
a $\Gamma$-Kotzig array of size $j\times k$, we simply take $j/2$ of $\Gamma$-Kotzig arrays of size $2\times k$ and "glue" them into an array $j\times k$.

Suppose $j$ is odd and there exists $\Gamma$-Kotzig array for a group $\Gamma$ containing only one involution $\iota$.
Recall that $kg=g_0$ for any $g\in \Gamma$.  On one hand $\sum_{i=1}^k(x_{1,i}+x_{2,i}+\ldots+x_{j,i})=g_0$, because  for some $\mu\in \Gamma$ we have $x_{1,i}+x_{2,i}+\ldots+x_{j,i}=\mu$ for any $i=1,2,\ldots,k$. On the other $\sum_{i=1}^k(x_{1,i}+x_{2,i}+\ldots+x_{j,i})=j\sum_{g\in \Gamma}g=k\iota=\iota$, a contradiction.
 Assume now that $\Gamma\in \gr$. To construct
a $\Gamma$-Kotzig array of size $j\times k$, we simply take $\lfloor j/2 \rfloor-1$ of $\Gamma$-Kotzig arrays of size $2\times k$, one $\Gamma$-Kotzig array of size $3\times k$ and "glue" them into an array $j\times k$.~\qed



\section{Final remarks}
Since any group $\Gamma\in\gr$ has $m$-ZSP for any $m\geq3$,  we state the following conjecture.
\begin{con}Let $\Gamma\in\gr$ have even order and $2^t-1\geq3$ involutions. For
every partition $n-1 = r_1 + r_2 + \ldots + r_t$ of $n-1$, with $r_i \geq 3$ for $1 \leq i \leq t$ and for any possible positive integer $t$, there is a partition of $\Gamma-\{g_0\}$ into pairwise disjoint subsets $A_1, A_2,\ldots , A_t$, such that $|A_i| = r_i$ and $\sum_{a\in A_i}a = g_0$ for $1 \leq i \leq t$. 
\end{con}

\end{document}